\documentclass[11pt,reqno]{amsart}
\usepackage{amssymb}
\usepackage{amsmath}
\usepackage{amsfonts}

\newtheorem{theorem}{Theorem}
\theoremstyle{plain}

\newtheorem{lemma}{Lemma}

\newtheorem{proposition}{Proposition}
\newtheorem{remark}{Remark}

\numberwithin{equation}{section}
\numberwithin{equation}{section}
\input{tcilatex}

\begin{document}
\title[Energy decay]{Stabilization of the wave equation with external force. 
}
\author{M. Daoulatli}
\address{Department of Mathematics and Informatics, ISSATS, University of
Sousse,\\
Sousse 4003 , Tunisia.}
\email[M. Daoulatli]{moez.daoulatli@infcom.rnu.tn}
\date{\today }

\begin{abstract}
We study the rate of decay of the energy functional of solutions of the wave
equation with localized damping and a external force. We prove that the
decay rates of the energy functional \ is determined from a forced
differential equation.
\end{abstract}

\maketitle

\section{Introduction}

This article is devoted to the study of stabilization for the wave equation
with external force on a compact Riemannian manifold with boundary. In the
first part of this paper, we consider the following wave equation with
linear internal damping and external force%
\begin{equation}
\left\{ 
\begin{array}{ll}
\partial _{t}^{2}u-\Delta u+a\left( x\right) \partial _{t}u=f\left(
t,x\right) & 
\mathbb{R}
_{+}\times M \\ 
u=0 & 
\mathbb{R}
_{+}\times \partial M \\ 
\left( u\left( 0\right) ,\partial _{t}u\left( 0\right) \right) =\left(
u_{0},u_{1}\right) & 
\end{array}%
\right.  \label{sys:linear}
\end{equation}%
Here $M=\left( M,q\right) $ is a compact, connected Riemannian manifold of
dimension $d$, with $C^{\infty }$ boundary $\partial M$, where $q$ denotes a
Riemannian metric of class $C^{\infty }$. $\Delta $ the Laplace--Beltrami
operator on $M.$ $a\left( x\right) $ is a non negative function in $%
C^{\infty }\left( M\right) $ and $f$ is a function in $L^{2}\left( 
\mathbb{R}
_{+}\times M\right) .$

We define the energy space%
\begin{equation*}
\mathcal{H}=H_{0}^{1}\left( M\right) \times L^{2}(M)
\end{equation*}%
where%
\begin{equation*}
H_{0}^{1}\left( M\right) =\left\{ u\in H^{1}\left( M\right) ;\text{ }%
u|_{\partial M}=0\right\}
\end{equation*}%
which is a Hilbert space. Linear semigroup theory applied to (\ref%
{sys:linear}), provides the existence of a unique solution $u$ in the class%
\begin{equation*}
u\in C^{0}\left( 
\mathbb{R}
_{+},H_{0}^{1}\left( M\right) \right) \cap C^{1}\left( 
\mathbb{R}
_{+},L^{2}\left( M\right) \right)
\end{equation*}%
With $\left( \text{\ref{sys:linear}}\right) $ we associate the energy
functional given by%
\begin{equation*}
E_{u}\left( t\right) =\frac{1}{2}\int_{M}\left\vert \nabla u\left(
t,x\right) \right\vert ^{2}+\left\vert \partial _{t}u\left( t,x\right)
\right\vert ^{2}dx.
\end{equation*}%
The energy $E_{u}(t)$ is topologically equivalent to the norm on the space $%
\mathcal{H}.$ Under these assumptions, the energy functional satisfies the
following identity%
\begin{equation}
E_{u}\left( t\right) +\int_{s}^{t}\int_{M}a\left( x\right) \left\vert
\partial _{t}u\right\vert ^{2}dxd\sigma =E_{u}\left( s\right)
+\int_{s}^{t}\int_{M}f\partial _{t}udxd\sigma  \label{energy identity}
\end{equation}%
for every $t\geq s\geq 0.$

The topic of interest is rate of decay of the energy functional. This
problem has a very long history. The connection between controllability,
observability and stabilization was discovered \cite{russell} and
effectively used in the context of linear PDE systems.

When $f=0$ and the damping term acts on the hole manifold the problem has
been studied by many authors \cite{andrade} and reference therein. For the
wave equation with localized linear damping term, we mention the works of
Rauch- Taylor\cite{rauch taylor} and Bardos et al \cite{blr} in which
microlocal techniques is used. In particular the notion of geometric
control. We cite also the works of Lasiecka et al \cite{las-tri-ya} and
Triggiani- Yao \cite{tri-yao} in which another approach based on Remannian
geometry is presented.

Particular attention has been paid to the case when $M$ is a bounded domain
and the damping is linearly bounded \cite{las1}\cite{komor} and reference
therein$.$ Under certain geometric condition, the energy functional decays
exponentially. Damping that does not satisfy such linear bound near the
origin (e.g. when the damping has polynomial, exponential or logarithmic
behavior near the origin) results in a weaker form of the energy that could
be expressed by algebraic, logarithmic (or possibly slower) rates \cite%
{las-tat}\cite{fab}\cite{mart}$.$ Finally we mention the work of cavalcanti
et al \cite{cavalcante} when $M$ is a compact manifold with or without
boundary.

When $f\neq 0,$ the literature is less furnished, we specially mention the
works of Haraux \cite{har} and Zhu \cite{zhu} when the damping is globally
distributed.

We should also remark when the support of the dissipation may be arbitrarily
small require more regular initial data and result in very slow (logarithmic
or slower ) decay rates as shown in \cite{daou1} and reference therein.

We assume that the geodesics of $\bar{M}$ have no contact of infinite order
with $\partial M.$ Let $\omega $ be an open subset of $M$ and consider the
following assumption:

(G) $\left( \omega ,T\right) $ geometrically controls $M$, i.e. every
generalized geodesic of $M$, travelling with speed $1$ and issued at $t=0$,
enters the set $\omega $ in a time $t<T$. \ 

This condition is called Geometric Control Condition (see e.g. \cite{blr}) \
We shall relate the open subset $\omega $ with the damper $a$ by 
\begin{equation*}
\omega =\{x\in M:a(x)>0\}.
\end{equation*}

Under the assumption (G) it was proved in \cite{blr,lebeau}, that the energy
decays exponentially, moreover if there exits a maximal generalized geodesic
of $M$ that never meets the support of the damper $a,$ then we don't have
the exponential decay of the energy for initial data in the energy space. 

It is known that the exponential decay of the energy is equivalent to the
following observability inequality:

\begin{description}
\item[(A) Linear Observability inequality] There exist positive constants $T$
and $\alpha =\alpha (T)$, such that for every initial condition $\varphi
=\left( u_{0},u_{1}\right) \in \mathcal{H}$ the corresponding solution\
satisfies 
\begin{equation}
E_{v}(t)\leq \alpha \int_{t}^{t+T}\int_{\Omega }a(x)\left\vert \partial
_{t}v\right\vert ^{2}dxds.
\end{equation}%
for every $t\geq 0.$
\end{description}

In this paper, under the assumption (G), we show that for the non autonomous
case the corresponding observability inequality reads as follows:

\begin{description}
\item[(\textbf{B}) Non autonomous linear Observability inequality] There
exist positive constants $T$ and $\alpha =\alpha (T)$, such that for every
initial condition $\varphi =\left( u_{0},u_{1}\right) \in \mathcal{H}$ the
corresponding solution\ satisfies 
\begin{equation}
E_{v}(t)\leq \alpha \int_{t}^{t+T}\int_{M}a(x)\left\vert \partial
_{t}v\right\vert ^{2}+\left\vert f\left( s,x\right) \right\vert ^{2}dxds.
\end{equation}%
for every $t\geq 0.$
\end{description}

From the observability inequality above, we infer that the rate of decay of
the energy will depends on $\int_{M}\left\vert f\left( t,x\right)
\right\vert ^{2}dx.$ Now we state the main result of the first part of the
paper:

\begin{theorem}
\label{t:1}Let $u(t)$ is the solution to the linear problem (\ref{sys:linear}%
) with initial condition $\left( u_{0},u_{1}\right) \in \mathcal{H}$. We
assume that $\left( \omega ,T\right) $ satisfies the assumption (G) and%
\begin{equation*}
\Gamma \left( t\right) =C_{1,T}\int_{M}\left\vert f\left( t,x\right)
\right\vert ^{2}dx\in L^{1}\left( 
\mathbb{R}
_{+}\right) 
\end{equation*}%
with $C_{1,T}\geq 1$. Then%
\begin{equation*}
E_{u}\left( t\right) \leq 4e^{T}\left( S\left( t-T\right)
+\int_{t-T}^{t}\Gamma \left( s\right) ds\right) ,\qquad t\geq T
\end{equation*}%
\ where $S\left( t\right) $ is the solution of the following ordinary
differential equation 
\begin{equation}
\frac{dS}{dt}+\frac{1}{TC_{T}}S=\Gamma \left( t\right) ,\qquad S\left(
0\right) =E_{u}\left( 0\right) .  \label{sharp ODE}
\end{equation}%
where $C_{T}\geq 1.$
\end{theorem}

\subsection{\textbf{Applications for the linear case}}

Setting 
\begin{equation*}
\Gamma \left( t\right) =C_{1,T}\int_{M}\left\vert f\left( t,x\right)
\right\vert ^{2}dx
\end{equation*}%
with $C_{1,T}\geq 1$.

The ODE $\left( \ref{sharp ODE}\right) $ governing the energy bound reduces
to%
\begin{equation}
\frac{dS}{dt}+CS=\Gamma \left( t\right)  \label{equation linear}
\end{equation}%
where constant $C>0$ does not depend on $E_{u}\left( 0\right) .$

\begin{enumerate}
\item If there are constants $M>0$ and $\theta >0,$ such that%
\begin{equation*}
\Gamma \left( t\right) \leq Me^{-\theta t}
\end{equation*}%
We have%
\begin{equation*}
\int_{t-T}^{t}e^{-\theta s}ds\leq \frac{1}{\theta }\left[ e^{\theta T}-1%
\right] e^{-\theta t},t\geq T
\end{equation*}%
Multiply (\ref{equation linear}) both sides by $\exp (Ct)$ and integrate
from $0$ to $t$, we obtain

\begin{enumerate}
\item $C>\theta $%
\begin{equation*}
E_{u}\left( t\right) \leq c\left( 1+E_{u}\left( 0\right) \right) e^{-\theta
t},t\geq 0
\end{equation*}

\item $C=\theta $%
\begin{equation*}
E_{u}\left( t\right) \leq c\left( 1+E_{u}\left( 0\right) \right) \left(
1+t\right) e^{-\theta t},t\geq 0
\end{equation*}

\item $C<\theta $%
\begin{equation*}
E_{u}\left( t\right) \leq c\left( 1+E_{u}\left( 0\right) \right)
e^{-Ct},t\geq 0
\end{equation*}
\end{enumerate}

\item If there are constants $M>0$ and $\theta >1,$ such that%
\begin{equation*}
\Gamma \left( t\right) \leq M\left( 1+t\right) ^{-\theta }
\end{equation*}
\ We have%
\begin{equation*}
\int_{t-T}^{t}\left( 1+s\right) ^{-\theta }ds\leq T\left( 1+t-T\right)
^{-\theta },t\geq T
\end{equation*}%
In order to obtain the rate of decay in this case, we use proposition \ref%
{lemma ode}. Then%
\begin{equation*}
E_{u}\left( t\right) \leq c\left( 1+t-T\right) ^{-\theta },t\geq T
\end{equation*}%
where $c>0$ and depends on $E_{u}\left( 0\right) .$
\end{enumerate}

\begin{remark}
If we consider the following system%
\begin{equation*}
\left\{ 
\begin{array}{ll}
\partial _{t}^{2}u-\Delta u+a\left( x\right) g\left( \partial _{t}u\right)
=f\left( t,x\right)  & 
\mathbb{R}
_{+}\times M \\ 
u=0 & 
\mathbb{R}
_{+}\times \partial M \\ 
\left( u\left( 0\right) ,\partial _{t}u\left( 0\right) \right) =\left(
u_{0},u_{1}\right) \in H_{0}^{1}\left( M\right) \times L^{2}\left( M\right) 
& 
\end{array}%
\right. 
\end{equation*}%
with $g$ continuous, monotone increasing function, vanishing at the origin
and linearly bounded. Then the result of the theorem above remains true.
\end{remark}

\subsection{The nonlinear case}

In the second part of the paper we study the rate of decay of the energy
functional of solution of the wave equation with nonlinear damping and
external force. More precisely, we consider the following system%
\begin{equation}
\left\{ 
\begin{array}{ll}
\partial _{t}^{2}u-\Delta u+a\left( x\right) g\left( \partial _{t}u\right)
=f\left( t,x\right) & 
\mathbb{R}
_{+}\times M \\ 
u=0 & 
\mathbb{R}
_{+}\times \partial M \\ 
\left( u\left( 0\right) ,\partial _{t}u\left( 0\right) \right) =\left(
u_{0},u_{1}\right) \in H_{0}^{1}\left( M\right) \times L^{2}\left( M\right)
& 
\end{array}%
\right.  \label{sys:nonlinear}
\end{equation}

$g$ is a continuous, monotone increasing function vanishing at the origin.
Moreover we assume that, there exists a positive constant $m,$ such that%
\begin{equation}
\frac{1}{m}\left\vert s\right\vert ^{2}\leq g\left( s\right) s\leq
m\left\vert s\right\vert ^{2},\text{ }\left\vert s\right\vert >\eta 
\label{behavior infinity}
\end{equation}%
for some $\eta >0.$ $a\left( x\right) $ is a non negative function in $%
C^{\infty }\left( M\right) $ and $f$ is in $L^{2}\left( 
\mathbb{R}
_{+}\times M\right) .$

Nonlinear semigroup theory applied to (\ref{sys:linear}), provides the
existence of a unique solution $u$ in the class%
\begin{equation*}
u\in C^{0}\left( 
\mathbb{R}
_{+},H_{0}^{1}\left( M\right) \right) \cap C^{1}\left( 
\mathbb{R}
_{+},L^{2}\left( M\right) \right)
\end{equation*}%
Under these assumptions on the behavior on the damping, the energy
functional satisfies the following identity%
\begin{equation}
E_{u}\left( t\right) +\int_{s}^{t}\int_{M}a\left( x\right) g\left( \partial
_{t}u\right) \partial _{t}udxd\sigma =E_{u}\left( s\right)
+\int_{s}^{t}\int_{M}f\partial _{t}udxd\sigma
\end{equation}%
for every $t\geq s\geq 0.$

It is well known, for the nonlinear problem without a external force the
corresponding observability inequality \cite{las-tat,MID}... reads as
follows:

\begin{description}
\item[(C) Nonlinear Observability Inequality] There exists a constant $T>0$
and a concave, continuous, monotone increasing function $h:\mathbb{R}%
_{+}\rightarrow \mathbb{R}_{+}$, $h(0)=0$ (possibly dependent on $T$ ) such
that the solution $u(t,x)$ to the nonlinear problem (\ref{sys:nonlinear})
with initial data $\varphi =\left( u_{0},u_{1}\right) $ and $f\equiv 0$
satisfies 
\begin{equation}
E_{u}\left( t\right) \leq h\left( \int_{t}^{t+T}\int_{\Omega }a(x)g(\partial
_{t}u)\partial _{t}u\;dxds\right) ,  \label{observability nonlinear}
\end{equation}%
for every $t\geq 0.$
\end{description}

The function $h(s)$ in (\ref{observability nonlinear}) depends on the
nonlinear map $g(s)$, and ultimately determines the decay rates for the
energy $E_{u}\left( t\right) $. The energy decay for the \emph{nonlinear}
problem will be determined from the following ODE 
\begin{equation}
S_{t}+h^{-1}\left( CS\right) =0,\quad S\left( 0\right) =E_{u}\left( 0\right)
\label{ODE}
\end{equation}%
we show that under the assumption (\textbf{G) }we obtain the following
observability inequality

\begin{description}
\item[(D) Nonlinear Non-autonomous Observability Inequality] There exists a
constant $T>0$ and a concave, continuous, monotone increasing function $h:%
\mathbb{R}_{+}\rightarrow \mathbb{R}_{+}$, $h(0)=0$ (possibly dependent on $%
T $ ) such that the solution $u(t,x)$ to the nonlinear problem (\ref%
{sys:nonlinear}) with initial data $\varphi =\left( u_{0},u_{1}\right) $
satisfies 
\begin{equation}
E_{u}\left( t\right) \leq h\left( \int_{t}^{t+T}\int_{\Omega }a(x)g(\partial
_{t}u)\partial _{t}u\;dxds+\int_{t}^{t+T}\int_{M}\left\vert f\left(
s,x\right) \right\vert ^{2}dxds\right) ,
\end{equation}%
for every $t\geq 0.$
\end{description}

Before giving the main result of this section, we will define some needed
functions. According to \cite{las-tat} there exists a strictly increasing
function $h_{0}$ with $h_{0}\left( 0\right) =0$ such that%
\begin{equation*}
h_{0}\left( g\left( s\right) s\right) \geq \epsilon _{0}\left( \left\vert
s\right\vert ^{2}+\left\vert g\left( s\right) \right\vert ^{2}\right) ,\text{
}\left\vert s\right\vert \leq \eta
\end{equation*}%
for some $\epsilon _{0},\eta >0.$ For the construction of such function we
refer the interested reader to \cite{las-tat,daou 2}. With this function, we
define%
\begin{equation*}
h=I+\mathfrak{m}_{a}\left( M_{T}\right) h_{0}\circ \frac{I}{\mathfrak{m}%
_{a}\left( M_{T}\right) }
\end{equation*}%
where 
\begin{equation*}
\mathfrak{m}_{a}=a\left( x\right) dxdt\text{ and }M_{T}=\left( 0,T\right)
\times M
\end{equation*}%
We can now proceed to state the main result of the second part of the paper

\begin{theorem}
\label{t:2}Let $u(t)$ is the solution to the nonlinear problem (\ref%
{sys:nonlinear}) with initial condition $\left( u_{0},u_{1}\right) \in 
\mathcal{H}$. We assume that $\left( \omega ,T\right) $ satisfies the
assumption (G) and 
\begin{equation*}
\Gamma \left( t\right) =2\int_{M}\left\vert f\left( t,x\right) \right\vert
^{2}dx+\psi ^{\ast }\left( \left\Vert f\left( t,.\right) \right\Vert
_{L^{2}\left( M\right) }\right) \in L^{1}\left( 
\mathbb{R}
_{+}\right) 
\end{equation*}%
where $\psi ^{\ast }$ is the convex conjugate of the function $\psi ,$
defined by%
\begin{equation*}
\psi \left( s\right) =\left\{ 
\begin{array}{lc}
\frac{1}{2T}h^{-1}\left( \frac{s^{2}}{8C_{T}e^{T}}\right)  & s\in 
\mathbb{R}
_{+} \\ 
+\infty  & s\in 
\mathbb{R}
_{-}^{\ast }%
\end{array}%
\right. 
\end{equation*}%
with $C_{T}\geq 1.$ Then%
\begin{equation*}
E_{u}\left( t\right) \leq 4e^{T}\left( S\left( t-T\right)
+\int_{t-T}^{t}\Gamma \left( s\right) ds\right) ,\qquad t\geq T
\end{equation*}%
where $S\left( t\right) $ is the solution of the following ordinary
differential equation 
\begin{equation}
\frac{dS}{dt}+\frac{1}{4T}h^{-1}\left( \frac{1}{K}S\right) =\Gamma \left(
t\right) ,\qquad S\left( 0\right) =E_{u}\left( 0\right) .
\end{equation}%
with, $K\geq C_{T}.$
\end{theorem}

\subsubsection{Applications for the nonlinear case}

\begin{proposition}
\label{lemma ode}Let $p$ a differentiable, positive strictly increasing
function on $%
\mathbb{R}
_{+},$ \ We assume that there exists $m_{1}>0$ such that, $p\left( x\right)
\leq m_{1}x$ for every $x\in \left[ 0,\eta \right] $ for some $0<\eta <<1$
and the following property 
\begin{equation}
p\left( Kx\right) \geq mp\left( K\right) p\left( y\right)
\label{Lem:p lower bound}
\end{equation}%
holds for some $m>0$ and for every $\left( K,x\right) \in \left[ 1,+\infty %
\right[ \times 
\mathbb{R}
_{+}.$ $\Gamma \in C^{1}\left( 
\mathbb{R}
_{+}\right) .$ Let $S$ satisfying the following differential inequality 
\begin{equation*}
\frac{dS}{dt}+p\left( S\right) \leq \Gamma \left( t\right) ,\text{ }S\left(
0\right) \geq 0.
\end{equation*}

\begin{enumerate}
\item $\Gamma \left( t\right) =0$ for every $t\geq 0.$ We assume that $%
S\left( 0\right) >0.$ Then 
\begin{equation*}
S\left( t\right) \leq \psi ^{-1}\left( t\right) ,\text{ for every }t\geq 0
\end{equation*}%
where%
\begin{equation*}
\psi \left( x\right) =\int_{x}^{S\left( 0\right) }\frac{ds}{p\left( s\right) 
}.
\end{equation*}%
$x\in \left] 0,S\left( 0\right) \right] .$

\item $\Gamma \left( t\right) >0$ for every $t\geq 0.$

\begin{enumerate}
\item There exist $c>0$ and $\kappa \geq 1$ such that%
\begin{equation}
\frac{d}{dt}p^{-1}\left( \Gamma \left( t\right) \right) +c\Gamma \left(
t\right) <0,\text{ for every }t\geq 0  \label{application lemma 1}
\end{equation}%
and%
\begin{equation}
\begin{array}{c}
mp\left( \kappa \right) -\kappa c-1\geq 0 \\ 
\kappa p^{-1}\circ \Gamma \left( 0\right) \geq S\left( 0\right)%
\end{array}
\label{application assumption}
\end{equation}%
Then%
\begin{equation*}
S\left( t\right) \leq \kappa \psi ^{-1}\left( ct\right) ,\text{ for every }%
t\geq 0
\end{equation*}%
where%
\begin{equation*}
\psi \left( x\right) =\int_{x}^{p^{-1}\circ \Gamma \left( 0\right) }\frac{ds%
}{p\left( s\right) }.
\end{equation*}%
$x\in \left] 0,p^{-1}\circ \Gamma \left( 0\right) \right] .$

\item There exist $c>0$ and $\kappa \geq 1$ such that%
\begin{equation*}
\frac{d}{dt}p^{-1}\left( \Gamma \left( t\right) \right) +c\Gamma \left(
t\right) \geq 0,\text{ for every }t\geq 0
\end{equation*}%
and 
\begin{equation*}
\begin{array}{c}
mp\left( \kappa \right) -c\kappa -1\geq 0 \\ 
\kappa p^{-1}\circ \Gamma \left( 0\right) \geq S\left( 0\right)%
\end{array}%
\end{equation*}%
Then%
\begin{equation*}
S\left( t\right) \leq \kappa p^{-1}\circ \Gamma \left( t\right) ,\text{ for
every }t\geq 0
\end{equation*}
\end{enumerate}
\end{enumerate}
\end{proposition}

\begin{proof}
$\left. {}\right. $

\begin{enumerate}
\item If $S\left( 0\right) =0,$ since $S$ is positive and decreasing then $%
S\left( t\right) =0$ for every $t\geq 0.$ We assume that $S\left( 0\right)
>0 $. Let $\psi ,$ the function defined by 
\begin{equation*}
\psi \left( x\right) =\int_{x}^{S\left( 0\right) }\frac{ds}{p\left( s\right) 
}.
\end{equation*}%
then $\psi $ is a strictly decreasing function on $\left( 0,S\left( 0\right)
\right) $ and $\underset{x\rightarrow 0}{\lim }\psi \left( x\right) =+\infty
.$ We have%
\begin{equation*}
\frac{d}{dt}\psi \circ S\left( t\right) \geq 1
\end{equation*}%
Integrating from $0$ to $t,$ we obtain%
\begin{equation*}
\psi \circ S\left( t\right) \geq t,t\geq 0
\end{equation*}%
since $\psi $ is decreasing%
\begin{equation*}
S\left( t\right) \leq \psi ^{-1}\left( t\right) ,t\geq 0
\end{equation*}

\item $\left. {}\right. $

\begin{enumerate}
\item Let $\psi ,$ the function defined by 
\begin{equation*}
\psi \left( x\right) =\int_{x}^{p^{-1}\circ \Gamma \left( 0\right) }\frac{ds%
}{p\left( s\right) }.
\end{equation*}%
then $\psi $ is a strictly decreasing function on $\left] 0,p^{-1}\circ
\Gamma \left( 0\right) \right] $ and $\underset{x\rightarrow 0}{\lim }\psi
\left( x\right) =+\infty .$ We have%
\begin{equation*}
\frac{d}{dt}\psi \circ p^{-1}\circ \Gamma \left( t\right) =-\frac{\frac{d}{dt%
}p^{-1}\left( \Gamma \left( t\right) \right) }{\Gamma \left( t\right) }
\end{equation*}%
from $\left( \ref{application lemma 1}\right) ,$ we infer that%
\begin{equation*}
\frac{d}{dt}\psi \circ p^{-1}\circ \Gamma \left( t\right) \geq c
\end{equation*}%
Integrating from $0$ to $t,$ we obtain%
\begin{equation*}
\psi \circ p^{-1}\circ \Gamma \left( t\right) \geq ct
\end{equation*}%
this gives%
\begin{equation}
\Gamma \left( t\right) \leq p\circ \psi ^{-1}\left( ct\right) ,\text{ for
every }t\geq 0  \label{application 11}
\end{equation}%
Setting%
\begin{equation*}
y\left( t\right) =\kappa \psi ^{-1}\left( ct\right) ,t\geq 0.
\end{equation*}%
We have%
\begin{equation*}
y^{\prime }\left( t\right) +p\left( y\left( t\right) \right) =-c\kappa
p\circ \psi ^{-1}\left( ct\right) +p\left( \kappa \left( \psi ^{-1}\left(
ct\right) \right) \right)
\end{equation*}%
Using $\left( \ref{Lem:p lower bound}\right) $ and $\left( \ref{application
11}\right) $ 
\begin{eqnarray*}
y^{\prime }\left( t\right) +p\left( y\left( t\right) \right) &\geq &\left(
mp\left( \kappa \right) -c\kappa \right) p\circ \psi ^{-1}\left( ct\right) \\
&\geq &\left( mp\left( \kappa \right) -c\kappa \right) \Gamma \left( t\right)
\end{eqnarray*}%
$\left( \ref{application assumption}\right) $ gives%
\begin{equation*}
\begin{array}{c}
y^{\prime }\left( t\right) +p\left( y\left( t\right) \right) \geq \Gamma
\left( t\right) \\ 
y\left( 0\right) \geq S\left( 0\right)%
\end{array}%
\end{equation*}%
the result follows from the following lemma

\begin{lemma}
\label{lemma application}Let $p_{i}$ $\left( i=1,2\right) $ a positive
strictly increasing function on $%
\mathbb{R}
_{+}.$ Suppose that $S$ and $y$ are absolutely continuous functions and
satisfy 
\begin{equation}
\frac{dS}{dt}+p_{1}\left( S\right) \leq \Gamma \left( t\right) \text{ on }%
\left[ 0,+\infty \right[ .  \label{ode lemma application 1}
\end{equation}%
and 
\begin{equation}
\frac{dy}{dt}+p_{2}\left( y\right) \geq \Gamma _{1}\left( t\right) \text{ on 
}\left[ 0,+\infty \right[ .  \label{ode lemma application 2}
\end{equation}%
where $\Gamma ;\Gamma _{1}$ $\in L^{1}([0;\infty ));$ $\Gamma _{1}\geq
\Gamma \geq 0;$ $p_{1}\geq p_{2}\geq 0$. In addition; if%
\begin{equation*}
y\left( 0\right) \geq S\left( 0\right)
\end{equation*}%
Then%
\begin{equation*}
y\left( t\right) \geq S\left( t\right) ;\text{ for }t\geq 0
\end{equation*}
\end{lemma}

First we finish the proof of the proposition, then we give the proof of the
lemma.

\item Setting 
\begin{equation*}
y\left( t\right) =\kappa p^{-1}\circ \Gamma \left( t\right) ,t\geq 0.
\end{equation*}%
We have 
\begin{equation*}
y^{\prime }\left( t\right) +p\left( y\left( t\right) \right) =\kappa \left(
p^{-1}\circ \Gamma \right) ^{\prime }\left( t\right) +p\left( \kappa
p^{-1}\circ \Gamma \left( t\right) \right)
\end{equation*}%
Using $\left( \ref{Lem:p lower bound}\right) $ and the fact that 
\begin{equation*}
\frac{d}{dt}p^{-1}\left( \Gamma \left( t\right) \right) +c\Gamma \left(
t\right) \geq 0
\end{equation*}%
for some $c>0,$ we obtain 
\begin{equation*}
y^{\prime }\left( t\right) +p\left( y\left( t\right) \right) \geq \left(
mp\left( \kappa \right) -c\kappa \right) \Gamma \left( t\right)
\end{equation*}%
$\left( \ref{application assumption}\right) $ gives%
\begin{equation*}
\begin{array}{c}
y^{\prime }\left( t\right) +p\left( y\left( t\right) \right) \geq \Gamma
\left( t\right) \\ 
y\left( 0\right) \geq S\left( 0\right)%
\end{array}%
\end{equation*}%
the result follows from lemma \ref{lemma application}
\end{enumerate}
\end{enumerate}
\end{proof}

The proof of lemma \ref{lemma application} is borrowed from \cite{zhu}

\begin{proof}[Proof of lemma \protect\ref{lemma application}]
Suppose that there exists $t_{0}$ in $\left[ 0,+\infty \right[ ,$ such that%
\begin{equation*}
S\left( t_{0}\right) =y\left( t_{0}\right) \text{ and }S\left( t\right)
>y\left( t\right) \text{ on }\left[ t_{0},t_{0}+\epsilon \right]
\end{equation*}%
for some $\epsilon >0.$ Integrate $\left( \ref{ode lemma application 1}%
\right) $ and $\left( \ref{ode lemma application 2}\right) $ from $t_{0}$ to 
$t_{0}+\epsilon ,$ we obtain%
\begin{equation*}
S\left( t_{0}+\epsilon \right) -S\left( t_{0}\right)
+\int_{t_{0}}^{t_{0}+\epsilon }p_{1}\left( S\left( t\right) \right) dt\leq
\int_{t_{0}}^{t_{0}+\epsilon }\Gamma \left( t\right) dt
\end{equation*}%
and%
\begin{equation*}
y\left( t_{0}+\epsilon \right) -y\left( t_{0}\right)
+\int_{t_{0}}^{t_{0}+\epsilon }p_{2}\left( y\left( t\right) \right) dt\geq
\int_{t_{0}}^{t_{0}+\epsilon }\Gamma _{1}\left( t\right) dt
\end{equation*}%
therefore%
\begin{equation*}
S\left( t_{0}+\epsilon \right) +\int_{t_{0}}^{t_{0}+\epsilon }p_{1}\left(
S\left( t\right) \right) dt\leq y\left( t_{0}+\epsilon \right)
+\int_{t_{0}}^{t_{0}+\epsilon }p_{2}\left( y\left( t\right) \right) dt
\end{equation*}%
which gives%
\begin{eqnarray*}
S\left( t_{0}+\epsilon \right) -y\left( t_{0}+\epsilon \right) &\leq
&\int_{t_{0}}^{t_{0}+\epsilon }p_{2}\left( y\left( t\right) \right)
-p_{1}\left( S\left( t\right) \right) dt \\
&\leq &\int_{t_{0}}^{t_{0}+\epsilon }p_{1}\left( y\left( t\right) \right)
-p_{1}\left( S\left( t\right) \right) dt\leq 0
\end{eqnarray*}%
which contradict the fact that $S\left( t\right) >y\left( t\right) $ on $%
\left[ t_{0},t_{0}+\epsilon \right] .$
\end{proof}

Setting 
\begin{equation*}
\Gamma \left( t\right) =2\int_{M}\left\vert f\left( t,x\right) \right\vert
^{2}dx+\psi ^{\ast }\left( \left\Vert f\left( t,.\right) \right\Vert
_{L^{2}\left( M\right) }\right) 
\end{equation*}%
where $\psi ^{\ast }$ is the convex conjugate of the function $\psi ,$
defined by%
\begin{equation*}
\psi \left( s\right) =\left\{ 
\begin{array}{lc}
\frac{1}{2T}h^{-1}\left( \frac{s^{2}}{8C_{T}e^{T}}\right)  & s\in 
\mathbb{R}
_{+} \\ 
+\infty  & s\in 
\mathbb{R}
_{-}^{\ast }%
\end{array}%
\right. 
\end{equation*}%
and 
\begin{equation*}
\psi ^{\ast }\left( s\right) =\underset{y\in 
\mathbb{R}
}{\sup }\left[ sy-\varphi \left( y\right) \right] 
\end{equation*}

\begin{description}
\item[Superlinear damping] Assume 
\begin{equation*}
g\left( s\right) =\left\{ 
\begin{array}{cc}
s^{2}e^{-\frac{1}{s^{2}}} & 0\leq s<1 \\ 
-s^{2}e^{-\frac{1}{s^{2}}} & -1<s<0%
\end{array}%
\right. .
\end{equation*}%
We choose $h_{0}^{-1}\left( s\right) =s^{3/2}e^{-\frac{1}{s}},$ $0<s<\eta <<1
$ and%
\begin{equation*}
K>>\max \left( E_{u}\left( 0\right) +\left\Vert \Gamma \right\Vert
_{L^{1}\left( 
\mathbb{R}
_{+}\right) },C_{T}\right) .
\end{equation*}%
We have%
\begin{equation*}
\psi ^{\ast }\left( \left\Vert f\left( t,.\right) \right\Vert _{L^{2}\left(
M\right) }\right) \leq C\left( \left\Vert f\left( t,.\right) \right\Vert
_{L^{2}\left( M\right) }\left\vert \ln \left( \left\Vert f\left( t,.\right)
\right\Vert _{L^{2}\left( M\right) }\right) \right\vert ^{-\frac{1}{2}%
}+\left\Vert f\left( t,.\right) \right\Vert _{L^{2}\left( M\right)
}^{2}\right) .
\end{equation*}%
The ODE $\left( \ref{sharp ODE}\right) $ governing the energy bound reduces
to%
\begin{equation*}
\frac{dS}{dt}+CS^{3/2}e^{-\frac{1}{S}}\leq \Gamma \left( t\right) 
\end{equation*}%
with $C>0$ depends on $E_{u}\left( 0\right) .$ If there are constants $M>0$
and $\theta >1,$ such that%
\begin{equation*}
\Gamma \left( t\right) \leq M\left( 1+t\right) ^{-\theta }
\end{equation*}%
then%
\begin{equation*}
E_{u}\left( t\right) \leq \frac{c_{0}}{\ln \left( ct+c_{1}\right) },t\geq T
\end{equation*}%
with $c,c_{0},c_{1}>0$. These constants may depend on $E_{u}\left( 0\right) .
$

\item[Sublinear near the origin] Assume $g\left( s\right) s\simeq \left\vert
s\right\vert ^{1+r_{0}},$ $\left\vert s\right\vert <1,$ $r_{0}\in \left(
0,1\right) .$ We choose $h_{0}\left( s\right) =s^{2r_{0}/\left(
1+r_{0}\right) }$ for $0\leq s\leq 1$ and%
\begin{equation*}
K>>\max \left( E_{u}\left( 0\right) +\left\Vert \Gamma \right\Vert
_{L^{1}\left( 
\mathbb{R}
_{+}\right) },C_{T}\right) .
\end{equation*}%
We have%
\begin{equation*}
\psi ^{\ast }\left( \left\Vert f\left( t,.\right) \right\Vert _{L^{2}\left(
M\right) }\right) \leq C\left( \left\Vert f\left( t,.\right) \right\Vert
_{L^{2}\left( M\right) }^{r_{0}+1}+\left\Vert f\left( t,.\right) \right\Vert
_{L^{2}\left( M\right) }^{2}\right) 
\end{equation*}%
The ODE $\left( \ref{sharp ODE}\right) $ governing the energy bound reduces
to%
\begin{equation*}
\frac{dS}{dt}+CS^{\left( 1+r_{0}\right) /2r_{0}}\leq \Gamma \left( t\right) 
\end{equation*}%
with $C>0$ depends on $E_{u}\left( 0\right) .$

\begin{enumerate}
\item If there are constants $M>0$ and $\theta >1,$ such that%
\begin{equation*}
\Gamma \left( t\right) \leq M\left( 1+t\right) ^{-\theta }
\end{equation*}%
Then

\begin{enumerate}
\item $\theta \in \left] 1,\frac{1+r_{0}}{1-r_{0}}\right] .$%
\begin{equation*}
E_{u}\left( t\right) \leq c\left( 1+t-T\right) ^{-\frac{2r_{0}\theta }{%
1+r_{0}}},t\geq T
\end{equation*}%
where $c>0.$

\item $\theta \geq \frac{1+r_{0}}{1-r_{0}}$%
\begin{equation*}
E_{u}\left( t\right) \leq c\left( t-T\right) ^{-\frac{2r_{0}}{1-r_{0}}},t>T
\end{equation*}%
with $c>0$ and depends on $E_{u}\left( 0\right) $.
\end{enumerate}

\item If there are constants $M>0$ and $\theta >0,$ such that%
\begin{equation*}
\Gamma \left( t\right) \leq Me^{-\theta t}
\end{equation*}%
Then%
\begin{equation*}
E_{u}\left( t\right) \leq c\left( t-T\right) ^{-\frac{2r_{0}}{1-r_{0}}},t>T
\end{equation*}
\end{enumerate}
\end{description}

with $c>0$ and depends on $E_{u}\left( 0\right) $

\section{The linear case: Proof of theorem \protect\ref{t:1}}

\subsection{Preliminary results}

\begin{proposition}
Let $u$ be a solution of $\left( \ref{sys:linear}\right) $ with initial data
in the energy space. Then%
\begin{equation}
E_{u}\left( t\right) \leq \left( 1+\frac{1}{\epsilon }\right) e^{\epsilon
\left( t-s\right) }\left( E_{u}\left( s\right)
+\int_{s}^{t}\int_{M}\left\vert f\left( \sigma ,x\right) \right\vert
^{2}dxd\sigma \right)   \label{energy bound linear}
\end{equation}%
for every $\epsilon >0$ and for every $t\geq s\geq 0.$
\end{proposition}

\begin{proof}
Let $t\geq s\geq 0.$ From the energy identity%
\begin{equation*}
E_{u}\left( t\right) \leq E_{u}\left( s\right)
+\int_{s}^{t}\int_{M}f\partial _{t}udxd\sigma 
\end{equation*}%
Using Young's inequality%
\begin{equation*}
E_{u}\left( t\right) \leq E_{u}\left( s\right) +\frac{1}{\epsilon }%
\int_{s}^{t}\int_{M}\left\vert f\left( \sigma ,x\right) \right\vert
^{2}dxd\sigma +\epsilon \int_{s}^{t}E_{u}\left( \sigma \right) d\sigma 
\end{equation*}%
for every $\epsilon >0.$ Now Gronwall's lemma, gives%
\begin{equation*}
E_{u}\left( t\right) \leq e^{\epsilon \left( t-s\right) }\left( E_{u}\left(
s\right) +\frac{1}{\epsilon }\int_{s}^{t}\int_{M}\left\vert f\left( \sigma
,x\right) \right\vert ^{2}dxd\sigma \right) .
\end{equation*}
\end{proof}

The result below is a generalisation of the comparison lemma of Lasiecka and
Tataru \cite{las-tat}.

\begin{lemma}
\label{lemma las tat}Let $T>0$ and

\begin{itemize}
\item $\Gamma \in L_{loc}^{1}\left( 
\mathbb{R}
_{+}\right) $ and, non negative. Setting $\delta \left( t\right)
=\int_{t}^{t+T}\Gamma \left( s\right) ds$.

\item $W\left( t\right) $ be a non negative, continuous function for $t\in 
\mathbb{R}_{+}$. Moreover we assume that there exists a positive, monotone,
increasing function $\alpha $ such that%
\begin{equation*}
W\left( t\right) \leq \alpha \left( t-s\right) \left[ W\left( s\right)
+\int_{s}^{t}\Gamma \left( \sigma \right) d\sigma \right] ,\text{ for every }%
t\geq s\geq 0.
\end{equation*}

\item Suppose that $\ell $ and $I-\ell :%
\mathbb{R}
_{+}\rightarrow \mathbb{R}$ are increasing functions with $\ell (0)=0$ and 
\begin{equation}
W\left( \left( m+1\right) T\right) +\ell \left\{ W\left( mT\right) +\delta
\left( mT\right) \right\} \leq W\left( mT\right) +\delta \left( mT\right) 
\label{lemma las tat inequality}
\end{equation}%
for $m=0,1,2,..$ where $\ell \left( s\right) $ does not depend on $m.$ Then%
\begin{equation*}
W\left( t\right) \leq \alpha \left( T\right) \left( S\left( t-T\right)
+2\int_{t-T}^{t}\Gamma \left( s\right) ds\right) ,\qquad \forall t\geq T
\end{equation*}%
where $S\left( t\right) $ is a positive solution of the following nonlinear
differential equation%
\begin{equation}
\frac{dS}{dt}+\frac{1}{T}\ell \left( S\right) =\Gamma \left( t\right) \text{ 
};\qquad S(0)=W(0).  \label{Ode lema}
\end{equation}
\end{itemize}
\end{lemma}

\begin{proof}
To prove this result we use induction. Assume that $W\left( mT\right) \leq
S\left( mT\right) $ and prove that $W\left( \left( m+1\right) T\right) \leq
S\left( \left( m+1\right) T\right) $ where $S\left( t\right) $ is the
solution of $\left( \text{\ref{Ode lema}}\right) .$

Integrating the equation $\left( \text{\ref{Ode lema}}\right) $ from $mT$ to 
$\left( m+1\right) T$ yields%
\begin{equation}
S\left( \left( m+1\right) T\right) =S\left( mT\right) -\frac{1}{T}%
\int_{mT}^{\left( m+1\right) T}\ell \left( S\left( t\right) \right)
dt+\delta \left( mT\right)   \label{lem las tat 1}
\end{equation}%
On the other hand, we have%
\begin{equation*}
\frac{d}{dt}\left( S-\int_{0}^{t}\Gamma \left( s\right) ds\right) =-\frac{1}{%
T}\ell \left( S\right) \leq 0.
\end{equation*}%
therefore, for $t_{1}\geq t_{2}$ 
\begin{equation*}
S\left( t_{1}\right) \leq S\left( t_{2}\right) +\int_{t_{2}}^{t_{1}}\Gamma
\left( s\right) ds
\end{equation*}%
the function $\ell $ is increasing 
\begin{eqnarray*}
\ell \left( S\left( t\right) \right)  &\leq &\ell \left( S\left( mt\right)
+\int_{mT}^{t}\Gamma \left( s\right) ds\right) ,\text{ for }mT\leq t\leq
\left( m+1\right) T \\
&\leq &\ell \left( S\left( mt\right) +\delta \left( mT\right) \right) 
\end{eqnarray*}%
Using now, $\left( \ref{lem las tat 1}\right) ,$ we obtain%
\begin{equation*}
S\left( \left( m+1\right) T\right) \geq S\left( mT\right) +\delta \left(
mT\right) -\ell \left( S\left( mt\right) +\delta \left( mT\right) \right) 
\end{equation*}%
Since the function, $I-\ell $ is increasing%
\begin{equation*}
S\left( \left( m+1\right) T\right) \geq W\left( mT\right) +\delta \left(
mT\right) -\ell \left( W\left( mt\right) +\delta \left( mT\right) \right) 
\end{equation*}%
$\left( \ref{lemma las tat inequality}\right) ,$ gives%
\begin{equation*}
S\left( \left( m+1\right) T\right) \geq W\left( \left( m+1\right) T\right) .
\end{equation*}%
Setting $t=mT+\tau ,$ with $0\leq \tau <T.$ Then we obtain%
\begin{eqnarray*}
W\left( t\right)  &\leq &\alpha \left( \tau \right) \left[ W\left( t-\tau
\right) +\int_{t-\tau }^{t}\Gamma \left( s\right) ds\right]  \\
&\leq &\alpha \left( \tau \right) \left( S\left( t-\tau \right)
+\int_{t-\tau }^{t}\Gamma \left( s\right) ds\right)  \\
&\leq &\alpha \left( T\right) \left( S\left( t-T\right)
+2\int_{t-T}^{t}\Gamma \left( s\right) ds\right) ,\text{ for every }t\geq T.
\end{eqnarray*}
\end{proof}

\begin{proposition}
We assume that $\left( \omega ,T\right) $ satisfies the assumption (G).\
Then there exists $\hat{C}_{T}>0,$ such that the following inequality%
\begin{equation}
E_{u}\left( t\right) \leq \hat{C}_{T}\left[ \int_{t}^{t+T}\int_{M}a\left(
x\right) \left\vert \partial _{t}u\right\vert ^{2}+\left\vert f\left(
s,x\right) \right\vert ^{2}dxds\right]   \label{observability linear}
\end{equation}%
holds for every $t\geq 0$, for every solution $u$ of $\left( \ref{sys:linear}%
\right) $ with initial data in the energy space $\mathcal{H}$, for every $f$
in $L^{2}\left( 
\mathbb{R}
_{+}\times M\right) $ $.$
\end{proposition}

\begin{proof}
To prove this result we argue by contradiction. We assume that there exist a
sequence $\left( u_{n}\right) _{n}$ solution of $\left( \text{\ref%
{sys:nonlinear}}\right) $ with initial data in the energy space, a
non-negative sequence $\left( t_{n}\right) _{n}$ and $f_{n}$\ in $%
L^{2}\left( 
\mathbb{R}
_{+}\times M\right) ,$ such that 
\begin{equation}
\begin{array}{l}
E_{u_{n}}\left( t_{n}\right) \geq n\int_{t_{n}}^{t_{n}+T}\int_{M}a\left(
x\right) \left\vert \partial _{t}u_{n}\right\vert ^{2}+\left\vert
f_{n}\left( t,x\right) \right\vert ^{2}dxdt,%
\end{array}
\label{contradiction argument}
\end{equation}%
Moreover, $u_{n}$ has the following regularity%
\begin{equation}
u_{n}\in C\left( 
\mathbb{R}
_{+},H_{0}^{1}\left( M\right) \right) \cap C^{1}\left( 
\mathbb{R}
_{+},L^{2}\left( M\right) \right) 
\end{equation}%
Setting $\alpha _{n}=\left( E_{u_{n}}\left( t_{n}\right) \right) ^{1/2}>0$
and $v_{n}\left( t,x\right) =\frac{u_{n}\left( t_{n}+t,x\right) }{\alpha _{n}%
}.$ Then $v_{n}$ satisfies%
\begin{equation}
\left\{ 
\begin{array}{ll}
\partial _{t}^{2}v_{n}-\Delta v_{n}+a\left( x\right) \partial _{t}v_{n}=%
\frac{1}{\alpha _{n}}f_{n}\left( t_{n}+t,x\right)  & 
\mathbb{R}
_{+}\times M \\ 
v_{n}=0 & 
\mathbb{R}
_{+}\times \partial M \\ 
\left( v_{n}\left( 0\right) ,\partial _{s}v_{n}\left( 0\right) \right) =%
\frac{1}{\alpha _{n}}\left( u_{n}\left( t_{n}\right) ,\partial
_{t}u_{n}\left( t_{n}\right) \right)  & 
\end{array}%
\right.   \label{sys:vn}
\end{equation}%
Moreover, 
\begin{equation*}
E_{v_{n}}\left( 0\right) =1
\end{equation*}%
and%
\begin{equation}
\begin{array}{l}
1\geq n\int_{0}^{T}\int_{M}a\left( x\right) \left\vert \partial
_{t}v_{n}\right\vert ^{2}+\left\vert \frac{1}{\alpha _{n}}f_{n}\left(
t_{n}+t,x\right) \right\vert ^{2}dxdt,%
\end{array}
\label{contradiction argument1}
\end{equation}%
From the inequality above,\ we infer that%
\begin{equation}
\begin{array}{l}
\int_{0}^{T}\int_{M}a\left( x\right) \left\vert \partial
_{t}v_{n}\right\vert ^{2}dxdt\underset{n\rightarrow \infty }{\rightarrow }0,
\\ 
\int_{0}^{T}\int_{M}\left\vert \frac{1}{\alpha _{n}}f_{n}\left(
t_{n}+t,x\right) \right\vert ^{2}dxdt\underset{n\rightarrow \infty }{%
\rightarrow }0%
\end{array}
\label{l2 contradiction argument}
\end{equation}%
We have%
\begin{equation}
v_{n}\in C\left( \left[ 0,T\right] ,H_{0}^{1}\left( M\right) \right) \cap
C^{1}\left( \left[ 0,T\right] ,L^{2}\left( M\right) \right) 
\end{equation}%
Therefore,%
\begin{equation}
E_{v_{n}}\left( t\right) =E_{v_{n}}\left( 0\right)
-\int_{0}^{T}\int_{M}a\left( x\right) \left\vert \partial
_{t}v_{n}\right\vert ^{2}dxdt+\int_{0}^{T}\int_{M}\frac{1}{\alpha _{n}}%
f_{n}\left( t_{n}+t,x\right) \partial _{t}v_{n}dxdt
\label{energy identity vn}
\end{equation}%
and using $\left( \text{\ref{energy bound linear}}\right) ,$ we infer that%
\begin{eqnarray*}
E_{v_{n}}\left( t\right)  &\leq &2e^{T}\left( E_{v_{n}}\left( 0\right)
+\int_{0}^{T}\int_{M}\left\vert \frac{1}{\alpha _{n}}f_{n}\left(
t_{n}+t,x\right) \right\vert ^{2}dxdt\right)  \\
E_{v_{n}}\left( t\right)  &\leq &2e^{T}\left( 1+\frac{1}{n}\right) ,\text{
for every }t\in \left[ 0,T\right] 
\end{eqnarray*}%
This estimate allows one to show that the sequence $\left( v_{n},\partial
_{t}v_{n}\right) $ is bounded in $L^{\infty }\left( \left( 0,T\right) ,%
\mathcal{H}\right) $ then it admits a subsequence still denoted by $\left(
v_{n},\partial _{t}v_{n}\right) $ that converges weakly-* to $\left(
v,\partial _{t}v\right) $ in $L^{\infty }\left( \left( 0,T\right) ,\mathcal{H%
}\right) .$ Passing to the limit in the system satisfied by $v_{n},$ we
obtain 
\begin{equation}
\left\{ 
\begin{array}{ll}
\partial _{t}^{2}v-\Delta v=0 & \left] 0,T\right[ \times M \\ 
\partial _{t}v=0 & \left] 0,T\right[ \times \omega  \\ 
\left( v\left( 0\right) ,\partial _{s}v\left( 0\right) \right) \in \mathcal{H%
} & 
\end{array}%
\right.   \label{sys limite}
\end{equation}%
and the solution $v$ is in the class 
\begin{equation*}
C\left( \left[ 0,T\right] ,H_{0}^{1}\left( M\right) \right) \cap C^{1}\left( %
\left[ 0,T\right] ,L^{2}\left( M\right) \right) 
\end{equation*}%
We deduce as in J. Rauch and M. Taylor \cite{rauch taylor} or C. Bardos, G.
Lebeau, J. Rauch \cite{blr} that the set of such solutions is finite
dimensional and admits an eigenvector $v$ for $\Delta $. By unique
continuation for second order elliptic operator, we get $\partial _{t}v=0$.
Multiplying the equation by $v$ and integrating, we obtain $v=0.$ Now we
prove that $v_{n}\rightarrow 0,$ in the strong topology of $%
H_{loc}^{1}\left( \left( 0,T\right) ,H^{1}\left( M\right) \right) .$ For
that we use the notion of microlocal defect measures. These measures were
introduced by P G\'{e}rard \cite{ge1} and L. Tartar \cite{tatar}. Let $\mu $
the microlocal defect measure associated to the sequence $\left(
v_{n}\right) .$ From $\left( \ref{l2 contradiction argument}\right) $ we
infer that the supprot of $\mu $ is contained in characteristic set of the
wave operator and it propagates along the geodesic flow (G. Lebeau \cite%
{lebeau}$)$. Therefore 
\begin{equation*}
v_{n}\rightarrow 0,\text{ }H_{loc}^{1}\left( \left( 0,T\right) ,H^{1}\left(
\omega \right) \right) 
\end{equation*}%
Now the assumption (G) combined with the propagation of $\mu $ along
geodesic flow, gives 
\begin{equation*}
v_{n}\rightarrow 0,\text{ }H_{loc}^{1}\left( \left( 0,T\right) ,H^{1}\left(
M\right) \right) 
\end{equation*}%
This gives%
\begin{equation}
\int_{0}^{T}\varphi \left( t\right) E_{v_{n}}\left( t\right) dt\underset{%
n\rightarrow \infty }{\rightarrow }0  \label{energy vn integral}
\end{equation}%
for every $\varphi \in C_{0}^{\infty }\left( \left[ 0,T\right] \right) .$ On
the other hand, $E_{v_{n}}\left( 0\right) =1,$ therefore, from $\left( \ref%
{energy identity vn}\right) $ and the fact that%
\begin{equation*}
\int_{0}^{T}\int_{M}\left\vert \partial _{t}v_{n}\right\vert ^{2}dxdt\leq
2Te^{T}\left( 1+\frac{1}{n}\right) 
\end{equation*}%
we deduce that $E_{v_{n}}\left( t\right) \underset{n\rightarrow \infty }{%
\rightarrow }1,$ for every $t\in \left[ 0,T\right] .$ Since $E_{v_{n}}\left(
t\right) \leq 2,$ by Lebesgue's dominated convergence theorem%
\begin{equation*}
\int_{0}^{T}\varphi \left( t\right) E_{v_{n}}\left( t\right) dt\underset{%
n\rightarrow \infty }{\rightarrow }\int_{0}^{T}\varphi \left( t\right) dt.
\end{equation*}%
for every $\varphi \in C_{0}^{\infty }\left( \left[ 0,T\right] \right) .$ We
obtain a contradiction by choosing $\varphi $ such that $\int_{0}^{T}\varphi
\left( t\right) dt>0.$
\end{proof}

\subsection{Proof of Theorem \protect\ref{t:1}.}

Let $u$ be a solution of $\left( \ref{sys:linear}\right) $ with initial data
in the energy space. From the energy identity we have%
\begin{equation*}
\int_{t}^{t+T}\int_{M}a\left( x\right) \left\vert \partial _{t}u\right\vert
^{2}dxdt=E_{u}\left( t\right) -E_{u}\left( t+T\right)
+\int_{t}^{t+T}\int_{M}f\left( s,x\right) \partial _{t}udxds
\end{equation*}%
therefore, using Young's inequality%
\begin{eqnarray*}
\int_{t}^{t+T}\int_{M}a\left( x\right) \left\vert \partial _{t}u\right\vert
^{2}dxdt &\leq &E_{u}\left( t\right) -E_{u}\left( t+T\right)  \\
&&+\epsilon \int_{t}^{t+T}\int_{M}\left\vert \partial _{t}u\right\vert
^{2}dxds+\frac{1}{\epsilon }\int_{t}^{t+T}\int_{M}\left\vert f\left(
s,x\right) \right\vert ^{2}dxds
\end{eqnarray*}%
for every $\epsilon >0.$ Now using the observability estimate $\left( \ref%
{observability linear}\right) $ and $\left( \ref{energy bound linear}\right)
,$ we can show that%
\begin{equation}
\int_{t}^{t+T}\int_{M}\left\vert \partial _{t}u\right\vert ^{2}dxds\leq
2Te^{T}\hat{C}_{T}\left[ \int_{t}^{t+T}\int_{M}a\left( x\right) \left\vert
\partial _{t}u\right\vert ^{2}+\left\vert f\left( s,x\right) \right\vert
^{2}dxds\right]   \label{observability}
\end{equation}%
Then setting $\epsilon =\frac{1}{4Te^{T}\hat{C}_{T}},$ we infer that%
\begin{equation*}
\frac{1}{2}\int_{t}^{t+T}\int_{M}a\left( x\right) \left\vert \partial
_{t}u\right\vert ^{2}dxdt\leq E_{u}\left( t\right) -E_{u}\left( t+T\right)
+\left( 1+Te^{T}\hat{C}_{T}\right) \int_{t}^{t+T}\int_{M}\left\vert f\left(
s,x\right) \right\vert ^{2}dxds
\end{equation*}%
Hence with $\tilde{C}_{T}=\left( 1+Te^{T}\hat{C}_{T}\right) ,$ we have%
\begin{equation*}
\int_{t}^{t+T}\int_{M}a\left( x\right) \left\vert \partial _{t}u\right\vert
^{2}dxdt\leq 2\left[ E_{u}\left( t\right) -E_{u}\left( t+T\right) +\tilde{C}%
_{T}\int_{t}^{t+T}\int_{M}\left\vert f\left( s,x\right) \right\vert ^{2}dxds%
\right] 
\end{equation*}%
Now from the observability estimate $\left( \ref{observability}\right) $%
\begin{equation*}
E_{u}\left( t\right) \leq 2\hat{C}_{T}\left[ E_{u}\left( t\right)
-E_{u}\left( t+T\right) +\left( \tilde{C}_{T}+1\right)
\int_{t}^{t+T}\int_{M}\left\vert f\left( s,x\right) \right\vert ^{2}dxds%
\right] 
\end{equation*}%
with $\hat{C}_{T}\geq 1.$ Setting $C_{1,T}=2\left( \tilde{C}_{T}+1\right) .$
We remark that%
\begin{equation*}
2\hat{C}_{T}\left[ E_{u}\left( t\right) -E_{u}\left( t+T\right)
+C_{1,T}\int_{t}^{t+T}\int_{M}\left\vert f\left( s,x\right) \right\vert
^{2}dxds\right] \geq C_{1,T}\int_{t}^{t+T}\int_{M}\left\vert f\left(
s,x\right) \right\vert ^{2}dxds
\end{equation*}%
Then for $C_{T}=4\hat{C}_{T}$, we have%
\begin{eqnarray*}
E_{u}\left( t\right) +C_{1,T}\int_{t}^{t+T}\int_{M}\left\vert f\left(
s,x\right) \right\vert ^{2}dxds &\leq &C_{T}\left[ E_{u}\left( t\right)
-E_{u}\left( t+T\right) \right.  \\
&&\left. C_{1,T}\int_{t}^{t+T}\int_{M}\left\vert f\left( s,x\right)
\right\vert ^{2}dxds\right] 
\end{eqnarray*}%
Therefore%
\begin{eqnarray*}
&&E_{u}\left( t+T\right) +\frac{1}{C_{T}}\left[ E_{u}\left( t\right)
+C_{1,T}\int_{t}^{t+T}\int_{M}\left\vert f\left( s,x\right) \right\vert
^{2}dxds\right]  \\
&\leq &E_{u}\left( t\right) +C_{1,T}\int_{t}^{t+T}\int_{M}\left\vert f\left(
s,x\right) \right\vert ^{2}dxds
\end{eqnarray*}%
Setting $t=mT,$ with $m\in 
\mathbb{N}
$ and using Lemma \ref{lem las tat 1}, we conclude that%
\begin{equation*}
E_{u}\left( t\right) \leq 4e^{T}\left( S\left( t-T\right)
+\int_{t-T}^{t}\Gamma \left( s\right) ds\right) ,\qquad \forall t\geq T
\end{equation*}%
where $S\left( t\right) $ is a positive solution of the following nonlinear
differential equation%
\begin{equation}
\frac{dS}{dt}+\frac{1}{C_{T}T}S=\Gamma \left( t\right) \text{ };\qquad
S(0)=E_{u}(0)
\end{equation}%
and 
\begin{equation*}
\Gamma \left( s\right) =C_{1,T}\int_{M}\left\vert f\left( s,x\right)
\right\vert ^{2}dx
\end{equation*}

\section{The nonlinear case: Proof of theorem \protect\ref{t:2}}

This part is devoted to the proof of theorem \ref{t:2}. First we give the
following energy inequality.

\begin{proposition}
Let $u$ be a solution of $\left( \ref{sys:nonlinear}\right) $ with initial
data in the energy space. Then the following inequality%
\begin{equation}
E_{u}\left( t\right) \leq \left( 1+\frac{1}{\epsilon }\right) e^{\epsilon
\left( t-s\right) }\left( E_{u}\left( s\right)
+\int_{s}^{t}\int_{M}\left\vert f\left( \sigma ,x\right) \right\vert
^{2}dxd\sigma \right)  \label{energy bound nonlinear}
\end{equation}%
holds for every $\epsilon >0$ and for every $t\geq s\geq 0.$
\end{proposition}

For the proof of $\left( \ref{energy bound nonlinear}\right) ,$ we have only
to proceed as in the proof of $\left( \ref{energy bound linear}\right) .$
Now we give the proof of theorem \ref{t:2}.

\begin{proof}[Proof of Theorem \protect\ref{t:2}]
Let $u$ be the solution of $\left( \ref{sys:nonlinear}\right) $ with initial
condition $\left( u_{0},u_{1}\right) $ in the energy space $\mathcal{H}.$
Let $t\geq 0$ and $\phi =u\left( t+\cdot \right) $ be the solution of 
\begin{equation}
\begin{array}{l}
\left\{ 
\begin{array}{ll}
\partial _{s}^{2}\phi -\Delta \phi +a\left( x\right) g\left( \partial
_{s}\phi \right) =f\left( t+s,x\right)  & 
\mathbb{R}
_{+}\times \Omega  \\ 
\phi =0 & 
\mathbb{R}
_{+}\times \partial \Omega  \\ 
\left( \phi \left( 0\right) ,\partial _{s}\phi \left( 0\right) \right)
=\left( u\left( t\right) ,\partial _{t}u\left( t\right) \right)  & 
\end{array}%
\right. 
\end{array}
\label{sys:nonautonomous translate}
\end{equation}%
We argue as in \cite{MID}. Define $z=\phi -v$, where $v$ is the solution of $%
\left( \ref{sys:linear}\right) $ with initial data $\left( u\left( t\right)
,\partial _{t}u\left( t\right) \right) $ and $f=f\left( t+\cdot ,\cdot
\right) .$\ Then $z$ satisfies the system 
\begin{equation*}
\left\{ 
\begin{array}{ll}
\partial _{t}^{2}z-\Delta z+a\left( x\right) g\left( \partial _{t}\phi
\right) -a(x)\partial _{t}v=0 & 
\mathbb{R}
_{+}\times \Omega  \\ 
z=0 & 
\mathbb{R}
_{+}\times \partial \Omega  \\ 
\left( z\left( 0\right) ,\partial _{t}z\left( 0\right) \right) =0 & 
\end{array}%
\right. 
\end{equation*}%
Let $T>0,$ such that $\left( \omega ,T\right) $ satisfies the assumption $%
\left( \text{G}\right) .$\ It is clear that $a\left( x\right) \left( g\left(
\partial _{t}\phi \right) -\partial _{t}v\right) \in L^{2}\left( \left(
0,T\right) \times \Omega \right) .$ This observation permits us to apply
energy identity, whence 
\begin{eqnarray*}
E_{z}(T) &=&\int_{0}^{T}\int_{M}a\left( x\right) \left( \partial
_{t}v-g(\partial _{t}\phi )\right) \partial _{t}z\;dxdt \\
&=&-\int_{0}^{T}\int_{M}a\left( x\right) \left( \left\vert \partial
_{t}v\right\vert ^{2}+g(\partial _{t}\phi )\partial _{t}\phi \right)
\;dxdt+\int_{0}^{T}\int_{M}g\left( \partial _{t}\phi \right) \partial
_{t}v+\partial _{t}v\partial _{t}\phi d\mathfrak{m}_{a}
\end{eqnarray*}%
The monotonicity of $g$ ($g\left( s\right) s\geq 0$) and the estimate above,
gives the following estimate: 
\begin{equation*}
\int_{0}^{T}\int_{M}a\left( x\right) \left( \left\vert \partial
_{t}v\right\vert ^{2}+g(\partial _{t}\phi )\partial _{t}\phi \right)
dxdt\leq \int_{0}^{T}\int_{M}g\left( \partial _{t}\phi \right) \partial
_{t}v+\partial _{t}v\partial _{t}\phi d\mathfrak{m}_{a}
\end{equation*}%
Now the observability estimate $\left( \ref{observability linear}\right) $,
gives 
\begin{equation*}
E_{u}\left( t\right) =E_{v}\left( 0\right) \leq \hat{C}_{T}\left(
\int_{0}^{T}\int_{M}g\left( \partial _{t}\phi \right) \partial
_{t}v+\partial _{t}v\partial _{t}\phi d\mathfrak{m}_{a}+\int_{0}^{T}\int_{M}%
\left\vert f\left( t+s,x\right) \right\vert ^{2}dxds\right) 
\end{equation*}%
for some $\hat{C}_{T}\geq 1.$ From the estimate above we infer that%
\begin{equation*}
E_{u}\left( t\right) \leq \hat{C}_{T}\left( \int_{t}^{t+T}\int_{M}g\left(
\partial _{t}u\right) \partial _{t}\tilde{v}+\partial _{t}\tilde{v}\partial
_{t}ud\mathfrak{m}_{a}+\int_{t}^{t+T}\int_{M}\left\vert f\left( s,x\right)
\right\vert ^{2}dxds\right) 
\end{equation*}%
where $\tilde{v}\left( s\right) =v\left( s-t\right) ,$ $s\geq t\geq 0.$

\begin{lemma}
\label{lem:concave estimates global}Setting%
\begin{equation*}
M_{s,t}=[s,t]\times \Omega ,\text{ }t\geq s\geq 0\text{ and }M_{0,t}=M_{t}
\end{equation*}%
Let $t\geq 0.$ For $i=0,1$ let 
\begin{equation*}
M_{t}^{0}=\left\{ \left( s,x\right) \in \left[ t,t+T\right] \times \Omega ;%
\text{ }\left\vert \partial _{s}u\left( s,x\right) \right\vert <\eta
_{0}\right\} ,\quad M_{t}^{1}=M_{t,t+T}\setminus M_{t}^{0}
\end{equation*}%
and define 
\begin{equation*}
\Theta \left( M_{t,t+t}\right) =\int_{M_{t}}\left\vert \partial _{s}u\left(
s\right) \partial _{s}v\left( s-t\right) \right\vert d\mathfrak{m}_{a},\quad
\Psi \left( M_{t}^{i}\right) =\int_{M_{t}^{i}}\left\vert g\left( \partial
_{s}u\left( s\right) \right) \partial _{s}v\left( s-t\right) \right\vert d%
\mathfrak{m}_{a},
\end{equation*}%
where $u$ and $v$ denotes respectively the solution of $\left( \text{\ref%
{sys:nonlinear}}\right) $ and $\left( \text{\ref{sys:linear}}\right) $ with
initial data $\left( u_{0},u_{1}\right) $ and $\left( \left(
v_{0},v_{1}\right) =u\left( t\right) ,\partial _{t}u\left( t\right) \right) $%
.

\begin{enumerate}
\item The following inequality holds for every $\epsilon >0,$%
\begin{eqnarray}
\medskip \Psi \left( M_{t}^{0}\right) +\Theta \left( M_{t}^{0}\right)  &\leq
&\epsilon E_{u}\left( t\right) +\frac{C}{\epsilon }\mathfrak{m}_{a}\left(
M_{T}\right) h_{0}\left( {\frac{1}{\mathfrak{m}_{a}\left( M_{T}\right) }}%
\int_{M_{t,t+T}}g\left( \partial _{s}u\right) \partial _{s}u\;d\mathfrak{m}%
_{a}\right)   \label{Omega zero estimate global} \\
&&+\epsilon \int_{t}^{t+T}\int_{M}\left\vert f\left( s,x\right) \right\vert
^{2}dxds
\end{eqnarray}%
with $C>0.$

\item Estimate on the damping near infinity. The following inequality%
\begin{eqnarray}
\medskip \Psi \left( M_{t}^{1}\right) +\Theta \left( M_{t}^{1}\right)  &\leq
&\epsilon E_{u}\left( t\right) +C\epsilon ^{-1}\left(
\int_{M_{t,t+T}}g\left( \partial _{s}u\right) \partial _{s}u\;d\mathfrak{m}%
_{a}\right)   \label{super nonlinear h1 global} \\
&&+\epsilon \int_{t}^{t+T}\int_{M}\left\vert f\left( s,x\right) \right\vert
^{2}dxds
\end{eqnarray}%
holds for every $\epsilon >0$ with $C>0$.
\end{enumerate}
\end{lemma}

For the proof of the lemma above we have only to proceed as in \cite[Lemma
3.1 cases 1 and 2]{MID} and to use the energy inequality $\left( \ref{energy
bound linear}\right) $. 

Now using $\left( \ref{Omega zero estimate global}\right) $ and $\left( \ref%
{super nonlinear h1 global}\right) ,$ we deduce that%
\begin{equation*}
E_{u}\left( t\right) \leq \tilde{C}_{T}\left( \epsilon E_{u}\left( t\right)
+C_{T,\epsilon }h\left( \int_{M_{t,t+T}}g\left( \partial _{s}u\right)
\partial _{s}ud\mathfrak{m}_{a}+\int_{M_{t,t+T}}\left\vert f\left(
s,x\right) \right\vert ^{2}dxds\right) \right) 
\end{equation*}%
for every $\epsilon >0,$ where the function $h=I+\mathfrak{m}_{a}\left(
M_{T}\right) h_{0}\circ {\frac{I}{\mathfrak{m}_{a}\left( M_{T}\right) }}$
and $\tilde{C}_{T}\geq 1$. Setting $\epsilon $ small enough, e.g. $\epsilon =%
\frac{1}{2\tilde{C}_{T}}$%
\begin{equation}
E_{u}\left( t\right) \leq C_{T}h\left( \int_{t}^{t+T}\int_{M}g\left(
\partial _{s}u\right) \partial _{s}ud\mathfrak{m}_{a}+\int_{t}^{t+T}\int_{M}%
\left\vert f\left( s,x\right) \right\vert ^{2}dxds\right)   \label{proof 1}
\end{equation}%
for some $C_{T}\geq 1.$ This gives%
\begin{equation}
E_{u}\left( t\right) +\int_{t}^{t+T}\int_{M}\left\vert f\left( s,x\right)
\right\vert ^{2}dxds\leq 2C_{T}h\left( \int_{t}^{t+T}\int_{M}g\left(
\partial _{s}u\right) \partial _{s}ud\mathfrak{m}_{a}+\int_{t}^{t+T}\int_{M}%
\left\vert f\left( s,x\right) \right\vert ^{2}dxds\right) 
\label{Proof: observability nonlinear}
\end{equation}%
On the other hand, the energy identity gives%
\begin{equation}
\int_{t}^{t+T}\int_{M}a\left( x\right) g\left( \partial _{t}u\right)
\partial _{t}udxd\sigma \leq E_{u}\left( t\right) -E_{u}\left( t+T\right)
+\int_{t}^{t+T}\int_{M}\left\vert f\left( \sigma ,x\right) \partial
_{t}u\right\vert dxd\sigma   \label{proof 2}
\end{equation}%
Let $\psi ,$ defined by%
\begin{equation*}
\psi \left( s\right) =\left\{ 
\begin{array}{lc}
\frac{1}{2T}h^{-1}\left( \frac{s^{2}}{8C_{T}e^{T}}\right)  & s\in 
\mathbb{R}
_{+} \\ 
+\infty  & s\in 
\mathbb{R}
_{-}^{\ast }%
\end{array}%
\right. 
\end{equation*}%
It is clear that $\psi $ convex and proper function. Hence, we can apply
Young's inequality \cite{rockfellar} 
\begin{eqnarray*}
\int_{t}^{t+T}\int_{M}\left\vert f\left( \sigma ,x\right) \partial
_{t}u\right\vert dxd\sigma  &\leq &\int_{t}^{t+T}\left\Vert f\left( \sigma
,.\right) \right\Vert _{L^{2}}\left\Vert \partial _{t}u\left( \sigma
,.\right) \right\Vert _{L^{2}}d\sigma  \\
&\leq &\int_{t}^{t+T}\psi ^{\ast }\left( \left\Vert f\left( \sigma ,.\right)
\right\Vert _{L^{2}}\right) +\psi \left( \left\Vert \partial _{t}u\left(
\sigma ,.\right) \right\Vert _{L^{2}}\right) d\sigma 
\end{eqnarray*}%
where $\psi ^{\ast }$ is the convex conjugate of the function $\psi ,$
defined by 
\begin{equation*}
\psi ^{\ast }\left( s\right) =\underset{y\in 
\mathbb{R}
}{\sup }\left[ sy-\varphi \left( y\right) \right] 
\end{equation*}%
Using the energy inequality $\left( \ref{energy bound nonlinear}\right) $
and the observability estimate $\left( \ref{Proof: observability nonlinear}%
\right) ,$ we infer that%
\begin{equation*}
\int_{t}^{t+T}\psi \left( \left\Vert \partial _{t}u\left( \sigma ,.\right)
\right\Vert _{L^{2}}\right) d\sigma \leq \frac{1}{2}\left(
\int_{t}^{t+T}\int_{M}g\left( \partial _{s}u\right) \partial _{s}ud\mathfrak{%
m}_{a}+\int_{t}^{t+T}\int_{M}\left\vert f\left( s,x\right) \right\vert
^{2}dxds\right) 
\end{equation*}%
\ then $\left( \ref{proof 2}\right) ,$ gives%
\begin{eqnarray}
\int_{t}^{t+T}\int_{M}a\left( x\right) g\left( \partial _{t}u\right)
\partial _{t}udxd\sigma  &\leq &2\left( E_{u}\left( t\right) -E_{u}\left(
t+T\right) +\int_{t}^{t+T}\int_{M}\left\vert f\left( s,x\right) \right\vert
^{2}dxds\right.  \\
&&\left. +\int_{t}^{t+T}\psi ^{\ast }\left( \left\Vert f\left( \sigma
,.\right) \right\Vert _{L^{2}}\right) d\sigma \right) 
\end{eqnarray}%
\ The inequality above combined with the observability estimate $\left( \ref%
{proof 1}\right) $ and the fact $h=I+\mathfrak{m}_{a}\left( M_{T}\right)
h_{0}\circ {\frac{I}{\mathfrak{m}_{a}\left( M_{T}\right) }}$ is icreasing$,$
gives%
\begin{equation*}
E_{u}\left( t\right) \leq C_{T}h\left( 4\left( E_{u}\left( t\right)
-E_{u}\left( t+T\right) +\int_{t}^{t+T}\int_{M}\left\vert f\left( s,x\right)
\right\vert ^{2}dxd\sigma +\int_{t}^{t+T}\psi ^{\ast }\left( \left\Vert
f\left( \sigma ,.\right) \right\Vert _{L^{2}}\right) d\sigma \right) \right) 
\end{equation*}%
Setting%
\begin{equation*}
\Gamma \left( s\right) =2\int_{M}\left\vert f\left( s,x\right) \right\vert
^{2}dx+\psi ^{\ast }\left( \left\Vert f\left( s,.\right) \right\Vert
_{L^{2}}\right) 
\end{equation*}%
Therefore%
\begin{equation*}
E_{u}\left( t\right) +\int_{t}^{t+T}\Gamma \left( s\right) ds\leq Kh\left(
4\left( E_{u}\left( t\right) -E_{u}\left( t+T\right) +\int_{t}^{t+T}\Gamma
\left( s\right) dxds\right) \right) 
\end{equation*}%
with $K\geq C_{T}.$ Setting $\theta \left( t\right) =$\ $\int_{t}^{t+T}%
\Gamma \left( s\right) ds.$ Thus%
\begin{equation*}
E_{u}\left( t+T\right) +\frac{1}{4}h^{-1}\left( \frac{1}{K}\left(
E_{u}\left( t\right) +\theta \left( t\right) \right) \right) \leq
E_{u}\left( t\right) +\theta \left( t\right) 
\end{equation*}%
for every $t\geq 0.$ Take $t=mt,$ $m\in 
\mathbb{N}
$%
\begin{equation*}
E_{u}\left( \left( m+1\right) T\right) +\frac{1}{4}h^{-1}\left( \frac{1}{K}%
\left( E_{u}\left( mT\right) +\theta \left( mT\right) \right) \right) \leq
E_{u}\left( mT\right) +\theta \left( mT\right) 
\end{equation*}%
Setting $W\left( t\right) =E_{u}\left( t\right) ,$ $\ell \left( s\right) =%
\frac{1}{4}h^{-1}\circ \frac{I}{K}$ and $\Gamma \left( t\right)
=2\int_{M}\left\vert f\left( s,x\right) \right\vert ^{2}dx+\psi ^{\ast
}\left( \left\Vert f\left( s,.\right) \right\Vert _{L^{2}}\right) .$ It is
clear that the functions $\ell $ and $I-\ell $ are increasing on the
positive axis and $\ell \left( 0\right) =0.$ The function $\Gamma \in
L_{loc}^{1}\left( 
\mathbb{R}
_{+}\right) $ and non negative on $%
\mathbb{R}
_{+}$. According to lemma \ref{lemma las tat} 
\begin{equation*}
E_{u}\left( t\right) \leq 4e^{T}\left( S\left( t-T\right)
+\int_{t-T}^{t}\Gamma \left( s\right) ds\right) ,\qquad \forall t\geq T
\end{equation*}%
where $S\left( t\right) $ is the solution of the following nonlinear
differential equation%
\begin{equation}
\frac{dS}{dt}+\frac{1}{T}\ell \left( S\right) =\Gamma \left( t\right) \text{ 
};\qquad S(0)=W(0).
\end{equation}
\end{proof}

\end{document}